\documentclass[11pt]{amsart}
\usepackage{geometry}                
\usepackage[parfill]{parskip}    
\usepackage{graphicx}
\usepackage{amssymb}
\usepackage{epstopdf}
\usepackage[all, 2cell]{xy}
\usepackage[section]{placeins}
\usepackage[greek.polutoniko, english]{babel}

\usepackage{hyperref}
\usepackage{amsmath}
\usepackage{mathtools}
\usepackage{footmisc}

\title{Vladimir Voevodsky on the concept of mathematical structure in his letter exchange with Andrei Rodin}

\author{Andrei Rodin}

\date{\today}

\begin{document}
\begin{abstract}
In 2016 Vladimir Voevodsky sent the author an email message where he explained his conception of mathematical structure using a historical example borrowed from the \emph{Commentary to the First Book of Euclid's Elements} by Proclus; this message was followed by a short exchange where Vladimir clarified his conception of structure. In this Chapter Voevodsky's historical example is explained in detail, and the relevance of Voevodsky's conception of mathematical structure in Homotopy Type theory is shown. The Chapter also discusses some related historical and philosophical issues risen by Vladimir Voevodsky in the same email exchange. This includes a comparison of Voevodsky's conception of mathematical structure and other conceptions of structure found in the current literature. The concluding part of this Chapter includes relevant fragments of the email exchange between Vladimir Voevodsky and the author.       
\end{abstract}
\maketitle

\renewcommand{\contentsname}{} 
\tableofcontents

\section{Introduction} \label{intro}
I first met Vladimir Voevodsky in 2012 in Ljubljana  during a workshop\footnote{
Workshop on Higher Dimensional Algebra, Categories, and Types (Ljubljana, June 20, 2012) organised as a satellite event to the Fourth Workshop on Formal Topology (Ljubljana, June 15-19, 2012)}. Since then and until Vladimir's premature death in 2017 we met and talked at several occasions (including my visiting him in the Institute of Advanced Studies, Princeton, in 2015) and had a continuing email correspondence about mathematics, its history, and its philosophy. With the kind permission of Vladimir's widow Nadia Shalaby, I publish here a fragment of this correspondence where Vladimir puts forward an original conception of \emph{mathematical structure} and illustrates it with a historical example taken from the \emph{Commentary to the First Book of Euclid's Elements} \cite{Proclus:1873}, \cite{Proclus:1970}. 

The rest of the Chapter is organised as follows. In \ref{RSV} I present, explain and discuss Vladimir's argument providing historical and mathematical backgrounds, which are necessary for its understanding. More specifically, in \ref{PE} the relevant fragments of Euclid's \emph{Elements} and of Proclus' \emph{Commentary} are quoted and then explained. In \label{RU} Proclus' argument is reconstructed in modern logical terms: first, using the standard classical logical machinery and then using some more specific constructive logical means. In \ref{HoTT} I describe some basic ideas of Homotopy Type theory (HoTT) and show the relevance of Vladimir's historical example in HoTT. In \ref{FW} some further historical and philosophical implications of Vladimir's conception of mathematical structure and of his historical example are discussed. Finally,  \ref{AVC} comprises fragments of two Vladimir's messages and my replies in my English translation from the original Russian provided with some bibliographical references.

\section{Relations and Structures according to Voevodsky}  \label{RSV}
In his first message \ref{VR1} Vladimir claims that Proclus in his commentary on Euclid's definition of plane angle \cite[p.153]{Euclid:1908} clearly distinguishes between \emph{relations} and \emph{structures}\cite[p.99]{Proclus:1970}. What Vladimir says about relations suggests that he understands this concept in the usual way as a non-monadic predicate. But what Vladimir says about structures is less usual and calls for explanation. Indeed, as far as the concept of mathematical structure is understood \`a la Bourbaki as a family of sets with certain relations defined on these sets, it is not immediately clear how structures and relations can be compared. There is more than one way of thinking about mathematical structures in ontological terms but all the usual varieties of Mathematical Structuralism \cite{Hellman:2001} share the assumption according to which structures and relations in mathematics work closely together: in order to build the former one needs the latter. Think of Bourbaki's ``great types of structure'': the order structures, the algebraic structures, and the topological structures  \cite[p. 226-227]{Bourbaki:1950}. All of these are defined in terms of sets, families of sets, and relations;  the axiomatic foundations of the underlying set theory involve the primitive binary relation of membership. Some mathematical structuralists argue that mathematical structures can and probably should be conceived of as abstracted away from their underlying sets, i.e., conceived of in terms of ``pure'' relations abstracted from their relata
\footnote{In the 1950 manifesto written by Jean-Dieudonn\'e and signed with the name of Bourbaki the author points to set-theoretic paradoxes and related foundational difficulties and then says: ``The difficulties did not disappear until the notion of set itself disappeared (and with it all the metaphysical pseudo-problems concerning mathematical ``beings'') in the light of the recent work on logical formalism. From this new point of view, mathematical structures become, properly speaking, the only ``objects'' of mathematics.'' \cite[p. 225-226]{Bourbaki:1950}. Since the author does not provide here any precise reference it is difficult to say which contemporary works in logic, in his view, allow one to dispense with the concept of set in the foundations of mathematics. But he clearly states such a desideratum  anyway. A candidate theory, which has been later proposed as an alternative structuralist foundation of mathematics allowing one to dispense with sets in foundations, namely the Category theory, gained  a sufficient maturity for being considered to that role only in the mid-1960s \cite{Lawvere:1964},  \cite{Lawvere:1966}
}. But how one can possibly conceive of a mathematical structure without using relations? 

In his second message  \ref{VR2} Vladimir specifies that a relation on his account is a special and the most simple kind of structure; the idea of building a basic world picture (i.e., an ontology) in terms of objects and relations he describes as a ``laughably simplified version'' of a mathematical ontology built with mathematical structures. In \ref{VR1} Vladimir asks me a related historical question: When and how it happened that conceiving of the universe in terms of relations prevailed over conceiving of the universe in terms of structures? This question also appears to be wholly in odds with the received historiography of mathematics and the general history of ideas. While the history of relationist reasoning in mathematics and science can be traced back to Leibniz and Aristotle, the structuralist way of reasoning, according to the received view  \cite{Reck&Schiemer:2020}, did not emerge before the second half of the 19th century making its early appearance in Felix Klein's \emph{Erlangen Program} of 1872  \cite{Klein:1872}. 

I what follows I explain what Vladimir understands by a structure using his historical example. A discussion on Vladimir's historical question is postponed until \ref{FL}.

\subsection{Proclus on Euclid's Definition of Plane Angle} \label{PE}

Euclid's definition of plane angle (Definition 8 of Book 1 of his \emph{Elements}, Def.1.8. for short) is as follows:

\begin{quote}
[A] plane angle is the inclination of the lines, when two lines in a plane meet one another, and are not laid down straight-on with respect to one another. \cite[p.153]{Euclid:1908}
\end{quote}
 
Before we discuss Proclus' commentary on this definition let me make several commentaries of my own. Def.1.8.  is the first among several angle-related definitions found in the \emph{Elements}. In the next definition 1.9. Euclid defines a \emph{rectilinear} plane angle (a special case of 1.8. where the two lines are straight). Euclid's concept of straight line differs from today's concept bearing the same name because he assumes that a straight line is always bounded by its endpoints (cf. Def. 1.3), i.e., it is what today we would call a straight segment. Thus given two different straight lines $AB$, $CD$ intersecting in a general position in point $O$, Euclid distinguishes four different angles $\angle AOC$, $\angle COB$, $\angle BOD$ and $\angle DOA$ where all pairs of straight lines forming a given angle are different ($OA, OC$ for $\angle AOC$, etc.). Notice, that $\angle AOB$ and $\angle COD$ which today we qualify as \emph{straight} angles are not qualified as angles by Euclid: he explicitly rules out this possibility in the Def.1.8. Euclid also rules out the idea that two straight lines sharing an endpoint like $OA$ and $OC$ form two different angles (where one angle complements the other to the full circle): what Euclid calls an angle always has a radian measure $< \pi$. Euclid does not provide this latter detail in his definitions explicitly but this is how he uses the angle concept in his following Propositions and their proofs. 

 Def.1.8. admits for cases where one or both of the two lines are not straight. The only example of curvilinear angle found in the \emph{Elements} is the so-called \emph{horn angle} between a circle (or more precisely a circumference of a circle) and its tangent in the Proposition 3.16. This angle is formed by a line that is bounded (a straight line) and another line that has no boundary. 

In Book 11 which belongs to the stereometric part of the \emph{Elements} one finds more angle-related definitions (Def.11.5,6,7) including the definition of \emph{solid} (i.e. 3-dimensional) angle (Def.11.11). One does not find among these definitions, however, the concept of two-dimensional non-plane angle like an angle formed by two intersecting great circles of a sphere. But in his commentary on Euclid's Def.1.8. Proclus feels free to go beyond the formal limits of this and other Euclid's definitions and discusses examples of geometrical angles which he borrows from other sources. This includes non-planar two-dimensional angles (see below),  examples of angles formed by a single line self-intersecting line (cissoid and hyppopede \cite[p.103]{Proclus:1970}) and some other.

Proclus' commentary on Def.1.8. \cite[p.98-104]{Proclus:1970} starts with his attempt to place the concept of angle into one of three categories: relation, quality, and quantity. These three categories make part of Aristotle's list of 10 categories given in the very beginning of his \emph{Organon}. Proclus considers pros and cons for each of these three options and eventually comes to conclusion (supported by the authority of his predecessor as the Head of the (Neoplatonic) Academia Syrianus of Alexandria who apparently earlier expressed the same opinion) that the task is impossible and the only reasonable option is to treat the concept of angle as multi-categorical assuming that it somehow combines aspects of relation, aspects of quality, and aspects of quantity. Saying that angle is a quantity Proclus means that it is a magnitude like a straight line that allows one to establish that certain angles are equal, and that certain other angles are bigger or smaller than some yet other angles. Saying that angle is a quality Proclus means that being an angle is a qualification like qualifications ``straight'' and ``curve'' applied to a line (so he points here to a conception of angle as a ``broken line''). As an example of geometrical quality Proclus also points to the concept of geometrical \emph{figure}. So saying that angle is a quality amounts to saying that it is a geometrical figure of sort. This intuitive notion of geometrical figure, however, does not fall under Euclid's definition of (plane) figure 1.15, which requires a figure to be contained within certain boundaries. Finally, saying that angle is a relation Proclus interprets the ``inclination'' referred to by Euclid in his Def.1.8. as a \emph{sui generis} relation between the two lines mentioned in the same definition. This case is explained in detail in what follows\footnote{Otherwise Proclus' commentary on Def.1.8. can be described as an attempt to classify known definitions of angle using the aforementioned three Aristotle's categories. Proclus tentatively classifies Euclid's definition of plane angle into the category of relations but, as we shall shortly see, he rejects the idea that the concept of angle can be so defined. Proclus also stresses the fact that Euclid's definition does not cover the case of plane angles formed by self-intersecting lines (since Def.1.8. explicitly mentions two lines). Proclus' arguments seem to imply that Euclid's Def.1.8. is invalid but he stops short from stating this conclusion. He concludes instead in a neutral manner by saying ``So much we had to say about Euclid's definition, in part interpreting and in part exposing difficulties in it.'' \cite[p.104]{Proclus:1970}.}.

A more detailed analysis of the whole of Proclus' commentary on Def.1.8. is out of place here. Let me quote only a short fragment of this commentary where Proclus presents his argument against the idea of treating angles as relations, which is relevant to Vladimir's remark:

\begin{quote}
[I]f the angle is an inclination [as in Euclid's Def. 1.8.] and in general belongs to the class of relations, it will follow that, when the inclination is one, there is one angle and not more. For if the angle is nothing other than a relation between lines or between planes, how could there be one relation but many angles? If you imagine a cone cut by a triangle from apex to base, you will see one inclination at the apex of the half-cone, that of the sides of the triangle, but two separate angles, one the angle on the plane of the triangle, the other on the mixed surface of the cone; and both of these angles are contained by the above-mentioned two lines. The relation of these lines, then, did not make the angle. \cite[p.99]{Proclus:1970} 
\end{quote}

\begin{figure}[h]
\centering
\includegraphics[scale=0.5]{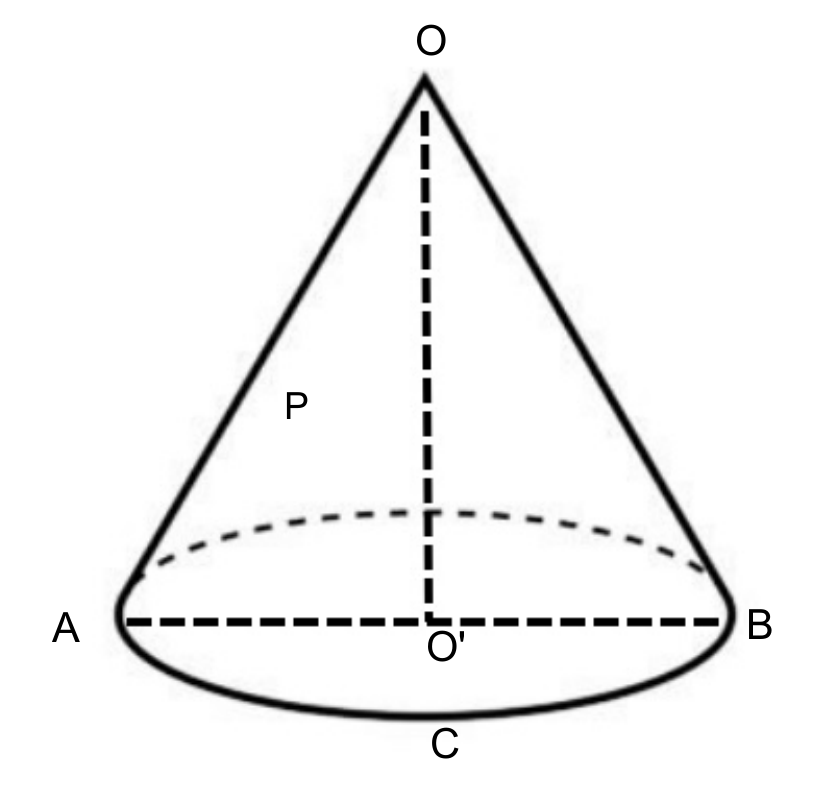}
\caption{Proclus' cone construction}
\label{fig:cone}
\end{figure}

In the language of today's school geometry the construction described by Proclus in the above quote can be presented as follows. Given cone $C$ with apex $O$ and base with centre $O'$ consider its triangular section $AOB$ by plane $P$ containing its axe $OO'$ where $AB$ is a diameter of the base. This determines plane angle $\angle AOB_{P}$ on $P$ that belongs to plane triangle $\triangle AOB_{P}$ and two non-plane angles $\angle AOB_{C}$, $\angle AOB'_{C}$  that belong to curve triangles $\triangle AOB_{C}$ and $\triangle AOB'_{C}$ formed on the cone surface (Fig. \ref{fig:cone}). 

Since Proclus considers only one of the two curve triangles we do the same in what follows. Proclus describes the obtained curved triangle $\triangle AOB_{C}$ as a``mixed surface'' meaning that two of its three sides are straight lines $AO$ and $BO$ but the third side is curve $\overset \frown {ACB}$, namely, a semicircle  (a half of the circumference of the base circle), cf. \cite[p.85]{Proclus:1970}. Notice that the above construction is obviously not available in the plane geometry where Euclid's Def.1.8. belongs. But this fact, in Proclus view,  does not make this construction irrelevant to his critique of Def.1.8. Since we want to reconstruct Proclus' argument charitably we may interpret his reasoning as follows: since (two-dimensional) angle cannot be defined as a relation of lines in the solid (3D) geometry  it should not be so defined in the plane geometry either because the plane geometry is a special case of the solid geometry (realised on each particular plane of \emph{the} Euclidean 3D space). 

Let us now see how Proclus' argument unfolds. A crucial element of this argument is the following assumption expressed by Proclus in the first two sentences of the above quote:  \emph{if} an angle is a pair of lines that hold the (binary) relation of \emph{inclination} \emph{then} given lines $OA, OB$ (that hold this relation) determine angle $\angle AOB$ \emph{uniquely}. Then Proclus applies his cone construction for showing that lines $OA, OB$, which presumingly do hold this relation, fail to determine a unique angle: we get instead two different angles $\angle AOB_{P}$  and $\angle AOB_{C}$. By contraposition Proclus concludes that angle cannot be defined (in the solid geometry and hence also in the plane geometry) as a pair of lines that hold the relation of inclination. 

It is clear that the above argument is schematic in the sense that it does not depend on what exactly is meant by the ``inclination'' of the two lines. The argument refutes the very idea that an angle can be defined as a pair of lines (i.e., the angle's sides) that stand in a certain binary relation to each other. The argument hinges on the aforementioned assumption, which should, in my view, be also understood  and reconstructed schematically as a general logical principle concerning binary relations. For further references I shall call this assumption the principle of \emph{Relational Uniqueness}, and denote it \textbf{RU} for short. In what follows two different modern reconstructions of  \ref{RU} are given: one in terms of Classical logic (see \ref{RUclass}) and the other in terms of Constructive logic (see \ref{RUconstr}) but without the full formalisation in either case. The classical reconstruction will play an auxiliary role but the constructive one will help us to see the analogy with HoTT, which Vladimir, to the best of my understanding, had in his mind when he pointed me to Proclus' commentary on Euclid's Def.1.8. in his message \ref{VR1}. In \ref{HoTT} this analogy  will be made explicit.

\subsection{Relational Uniqueness}  \label{RU}
Let me begin my reconstruction of \textbf{RU} with a historical remark. In his commentary on Def.1.8. Proclus applies Aristotle's concept of relation (\greektext{t\`a pr\'os ti}) \latintext which in many crucial respects differs from the modern notion \cite{Brower2016-BROAVC}. When Proclus applies the Aristotle's concept  he interprets it in his proper way. Having acknowledged this hermeneutic complexity, I do not engage myself in the present Chapter in a serious exegetic effort aiming at the historical reading of Proclus' \emph{Commentary} but limit my task to explaining and elaborating on the overtly anachronistic reading of Proclus suggested by Vladimir. I'll come back to the question of historical relevance of this reading in \ref{FL}.
The only remark concerning Aristotle's term (\greektext{t\`a pr\'os ti})\latintext, which is necessary in the present context, is this. As the term clearly suggests, Aristotle puts into the given category  ``things related'' (to some other things) rather than abstract relations disjoined from their relata. Accordingly, when Proclus reads Euclid's Def.1.8. as a relational definition of plane angle, he thinks of angle as a pair of lines that stand in the binary relation of \emph{inclination} (on an equal footing with a pair of parallel straight lines) rather than of the relation of inclination itself (construed as a binary predicate or otherwise). In other words, considering the option that ``angle is a relation'' Proclus, following Aristotle, considers the option that an angle is an \emph{instance} of a relation (between its sides). For further references I denote Euclid's relation of inclination referred to in his Def.1.8. by symbol  $\bowtie$ and write  $a \bowtie b$ when two lines $a,b$ hold this relation\footnote{What Euclid and Proclus say about the relation of inclination $\bowtie$ allows us to assume that this relation is symmetric.}.

Before I proceed to modern reconstructions of \textbf{RU} let me first tell what this principle is \underline{not}. When one defines a mathematical concept, say, that of geometrical \emph{square}, one does not, generally, make any assumption concerning the extension of the defined concept. Meaningful mathematical concepts  may have extensions containing a single object (like the concept of the empty set in the Set theory), finite or infinite extensions (like the concept of square in the elementary geometry) or empty extensions (like the concept of the largest prime number in arithmetic\footnote{This concept is certainly meaningful because the theorem that proves that the largest prime does not exist (see \emph{Elements} Prop. 9.20) is non-trivial and qualifies as a valuable piece of mathematical knowledge.}). Euclid's Def.1.8. does not aim at singling out an individual geometrical object, of course: there exist many different angles on the Euclidean plane. Thus \textbf{RU} does not concern the extension of Euclid's concept of plane angle defined with Def.1.8., which is obviously infinite. The Relational Uniqueness is a more subtle principle as we shall now see.


\subsubsection{Classical \textbf{RU}} \label{RUclass}

 \textbf{RU} has to do with how geometrical objects and constructions are identified (or ``individuated'' if one prefers) when they involve some relations. The issue identity of elementary geometrical objects is generally problematic. Think about the case of ``coinciding'' points for example. When two points coincide, are they the same point or two different points? Is the coincidence relation between points the same relation as the identity relation \cite{Rodin:2007}? We don't need to answer these difficult questions, however, in order to make sense of \textbf{RU}. For this limited purpose it is sufficient to assume that, given a relation, the identity conditions of its  \emph{relata} are already given; then \textbf{RU} tells us that the identity conditions of a composite object built from these relata are also specified. In other words, \textbf{RU} postulates that in order to specify the identity conditions for such a composite object it is sufficient to know the identity conditions of its components. 
 Consider parallel straight lines $a, b$, in symbols  $a \parallel b$\footnote{Here $a$ and $b$ are individual constants, i.e., proper names of certain well-individuated lines, but not variables.}. In that case \textbf{RU} implies that there exists unique composite object that can be described as \emph{pair} $(a,b)_{\parallel}$ of parallel lines. As shows Proclus using his cone construction, 2-dimensional angles (living in the 3D Euclidean space) do not share this property with parallels because given two lines $a,b$ such that $a \bowtie b$ one can get more than just one angle.

Thus we may tentatively formulate \textbf{RU} as follows: 

\begin{equation}
\text{Given objects $a, b$ such that $aRb$ there exists \textbf{unique} composite object (construction) $(a,b)_{R}$}.
\label{eq1}
\end{equation}

where $aRb$ is \emph{proposition} that says that $a$ and $b$ hold relation $R$. Notice that in order to be compatible with non-symmetric relations pair $(a,b)$ in \eqref{eq1} needs to be ordered
\footnote{If relation $R$ is non-symmetric then the order of its relata is fixed in the antecedent part of \eqref{eq1}, and thus the uniqueness of $(a,b)_{R}$ is still guaranteed. For generality we assume here that pair $(a,b)$ in \eqref{eq1} is ordered in all cases. If relation $R$ in  \eqref{eq1} is symmetric (like geometrical relations $\parallel$ and $\bowtie$) then $(a,b)_{R} = (b,a)_{R}$, i.e., the two composite objects are the same \label{fn1}}.  

Observe that since the identity conditions for $a$ and $b$ are given, one is in a position to define the identity conditions for the ordered pair $(a,b)$ in the obvious way:

\begin{equation}
(a,b) = (a',b')\ \text{if and only if}\ a=a'\ \text{and}\ b=b'
\label{eq2}
\end{equation}

This shows that the reference to  relation $R$ in  \eqref{eq1} is wholly redundant, so we get a more general principle:
 
\begin{equation}
\text{Given objects $a, b$ there exists \textbf{unique} composite object (ordered pair) $(a,b)$}.
\label{eq3}
\end{equation}

which can be justified in ZFC or another set theory using a standard representation of ordered pairs as sets. Thus the fact that given objects $a,b$ hold relation $R$ plays no role in how the ``composite object'' $(a,b)_{R}$ built of these given objects is identified. Under this reconstruction of \textbf{RU}, object  $(a,b)_{R}$ in  \eqref{eq1} is just the same as pair $(a,b)$, so the subscript $R$ in $(a,b)_{R}$ can be safely dropped. I call the reconstruction of \textbf{RU} via  \eqref{eq1} - \eqref{eq3}  \emph{classical} because it relies upon the Classical Predicate Calculus (using both its syntax and its standard Tarski-style semantics). 

 In view of \eqref{eq3} Proclus' argument can be paraphrased more formally as follows. Euclid's Def.1.8. says (according to Proclus' tentative interpretation) that angle $\widehat{a, b}$ is a construction that comprises lines $a,b$ such that  $a \bowtie b$, where $\bowtie$ is a binary relation. The semantics of binary relations involves the concept of ordered pair $(a,b)$. By \eqref{eq3}, given lines $a,b$, pair $(a,b)$ is unique\footref{fn1}. Thus given lines $a,b$, there can exist at most one angle of form $\widehat{a, b}$. But the cone construction shown at Fig.\ref{fig:cone} demonstrates that, given lines $a,b$, one can have two different angles of this form. Hence Euclid's Def.1.8. interpreted in terms of binary relation $\bowtie$ is not consistent with the common intuitive reasoning about geometrical angles. Hence  Def.1.8. (so interpreted) does not define the angle concept properly.  
  
 One can observe that the above classical reconstruction of \textbf{RU} trivialises this thesis in the sense that it doesn't suggest any possible epistemic scenario where it may fail. Further, this reconstruction does not explain the relevance of the concept of relation in \textbf{RU}: whether given objects $a,b$ hold some relation or not, there is unique (ordered) pair $(a,b)$. In these respects the classical reconstruction of  \textbf{RU} just given is hardly satisfactory (even if it perfectly saves Proclus' argument). In the next Section \ref{RUconstr} an alternative reconstruction of \textbf{RU} is proposed, which arguably represents Proclus' argument more faithfully. In addition, this alternative reconstruction of \textbf{RU} will bring us closer to HoTT and allow us to explain Vladimir's distinction between relations and structures that he attributes to Proclus in his message \ref{VR1}. 
 
Having said that, I would like to stress the relevance of the above classical reconstruction of Proclus' argument in the present discussion. It is relevant because it represents the common way of reasoning about relations in today's logic and mathematics. Vladimir's distinction between relations and structures involves a different view on relations, as we shall now see.


\subsubsection{Constructive \textbf{RU}} \label{RUconstr}
Before I propose a constructive reading of \textbf{RU} and of the rest of Proclus' argument let me present here some general reflections concerning a constructive treatment of relations, which motivate this approach. 
Consider a set of cities and binary relation ``to be train connected'' defined on that set. The fact that two cities, say, Paris and London, are train connected has a material evidence, namely, the Eurostar train line between the two cities. Thus the relation of train connectedness is not just an abstract matter: all train connected cities are materially related via a train line. Some other relations apparently are not materialised in a similar manner. For example, the fact that the Eiffel Tower is \emph{higher than} the Great Pyramid of Giza does not involve any obvious material link between the two buildings, which have been built independently in different times and in different places. The classical semantics of relations construed as binary predicates does not distinguish between such cases but treats binary relations uniformly \emph{sub specie aeternitatis}, i.e., from the viewpoint of the omniscient God, who always knows whether a given proposition is true or false without relying on any evidence for it.  

The irrelevance of epistemic matters in the Classical logic (including Classical Predicate Calculus) is an important ingredient of its usual philosophical underpinning (or motivation, if one prefers), which has been promoted by Frege \cite{Frege:1918} and his followers. This system of logic has been designed as a system of rules for reasoning about things as they supposedly \emph{are} without taking into account how these things are possibly \emph{known} to us. Taking this philosophical stance one is in a position to defend the Law of Excluded Middle (LEM) as a universal logical principle along the following line: quite independently of what we may know and what we may not know about the presence of Life on Mars, we are in a position to claim that either there is Life on Mars or there is no Life on Mars (provided that we understand unequivocally what is Life, what is ``being on Mars'', etc.). 

An alternative philosophical approach, which underpins various systems of  \emph{constructive} logic including Martin-L\"of's Type theory (MLTT)\cite{Martin-Lof:1984} that provides syntax and some crucial ingredients of informal semantics \cite{Martin-Lof:1987} for the standard version of HoTT \cite{UF:2013}, rejects the idea according to which one can meaningfully talk about true and false propositions without taking into account how these propositions are known. In the constructive logic\footnote{
I refer here to ``the constructive logic'' as a broad family of logical calculi motivated by epistemological ideas outlined in the main text. There are many different senses of being constructive, and there exist many different logical calculi that specify and implement the ``constructive approach'', broadly conceived, in many different ways. Among such logical calculi MLTT is particularly relevant to the present discussion because of its role in HoTT.} the classical notion of truth is replaced by the constructive notion of truth as an \emph{evidence} (aka witness, aka proof, aka truth-maker) \cite{Martin-Lof:1987}. Since many propositions (like ``there is Life on Mars''), to date, have neither a commonly accepted proof nor disproof, LEM does not hold in the constructive logic. 

Unlike the classical concept of truth the constructive concept of evidence is proposition-specific: while in the classical logic all propositions are evaluated with the same pair of truth values, evidence $e$ that proves proposition $P$, in symbols $e : P$, typically does not prove any other proposition. In MLTT this feature of constructive logic is syntactically realised via \emph{typing}: here term $e$ of type $P$, in symbols $e : P$, cannot be possibly a term of another type; a cross-identification of terms across different types is not allowed in MLTT.  Falsity and negation are more problematic matters in the constructive logic but the general idea is that the falsification of a given proposition is a (constructive) proof that this proposition has no proof.

Thus in order to assert that proposition $P$ is true \emph{constructively}, i.e., in order to make a \emph{judgement}, one needs to provide $P$ with some evidence (proof) $e$; the judgement in this case has logical form $e : P$.  In the above railroad example the relevant evidence supporting the claim that Paris and London are train connected, is the Eurostar line between the two cities. The Eiffel Tower example also admits for a constructive interpretation. Even if we assume that the Eiffel Tower (ET) and the Great Pyramid (GP) do not need to interact physically in order to sustain their independent existence, one can argue that the measurement and the comparison of the heights of the two buildings does require establishing such an interaction\footnote{More precisely, measuring and comparing the two buildings involves at least one more object, namely a \emph{knower} equipped with a measuring device and able to proceed the obtained empirical data appropriately. A simplified model of such a measurement involves a measuring stick that the knower applies to both buildings travelling from Paris to Giza. The causal networks of physical events that connects the knower to the two buildings via the measurement and the following computation may serve them as an evidence that ET is in fact higher than GP. In the social reality this causal network is by far more complex not only because it may involve more advanced measuring and computing tools but also because the knower is typically not a single person but rather a large network of different people living in different places at different times who may communicate with each other not only in the real time (and in the real place) but also remotely via written texts and other media. I must confess that I did not measure the heights of ET and GP myself, and I never even visited GP. Instead, I checked and compared the heights of the two buildings using Wikepedia. These details concerning the social machinery of human knowledge do not, however, constitute an objection to the epistemological thesis according to which in order to justify the claim that ET is higher than GP one needs to provide a conclusive evidence supporting this claim.}

Let us now return to geometrical examples. Let $a,b$ be straight lines,  and $a\otimes b$ be read  ``line $a$ intersects line $b$''. Point $p$ where the two lines intersect can be an evidence that the lines do intersect. Given such a point one is in a position to form judgement $p : a\otimes b$, i.e., to judge that the two lines intersect. The case of parallel lines is more difficult because in order to prove that two given lines $a,b$ are parallel it is necessary (and sufficient) to show that they 
 do \underline{not} intersect, i.e., that there is no $p$ such that $p : a\otimes b$. In this case an evidence cannot be a simple object like point $p$. An admissible evidence can consist, for example, of third straight line $c$ intersecting both $a$ and $b$, a proof that $c$ makes with  $a$ and $b$ equal corresponding angles $\alpha = \beta$, and, finally, a proof (using the assumption $\alpha = \beta$) that the existence of an intersection point of lines  $a$ and $b$ implies a contradiction\footnote{This follows (via a counterposition) from Euclid's Prop.1.17 that says that the sum of two angles of a triangle is strictly less then two right angles. Notice that Euclid's proof of Prop.1.17 does not use his Fifth Postulate, so this proposition is a theorem of the ``absolute'' elementary geometry, which constructively proves the existence of parallels in this theory. For a formal constructive treatment of elementary geometrical theories using the Type theory see \cite{vonPlato:1995}}).   

What can be an evidence that two lines $a, b$ meet the condition of Euclid's Def. 1.8., i.e., that they hold the relation of \emph{inclination} ($a\bowtie b$) referred to in this definition? A possible candidate is the angle $\widehat{a, b}$ itself, in symbols $\widehat{a, b} : a\bowtie b$. Of course, such a reading of Def.1.8. makes this definition circular. But this is not of our present concern \footnote{Many Euclid's definitions are similarly circular and can be read in both directions. Cf. for example Def.1.13: ``A boundary is that which is an extremity of anything''  \cite[p.153]{Euclid:1908}. Notice also that in Def.11.5 and Def.11.6 Euclid defines two more specific kinds of \emph{inclination} (one of a straight line to a plane and the other of a plane to another plane) in terms of plane angles defined with Def.1.8 and Def.1.9. For an overview of Euclid's definitions and their comparison with modern mathematical definitions see \cite[p. 39-41]{Mueller:1981}  }. Our (and Proclus') concern is whether or not the concept of geometrical angle has a content, which is not comprised in judgements of the form $\widehat{a, b} : a\bowtie b$. Using his cone construction (Fig.\ref{fig:cone}) Proclus shows that this is indeed the case. 

In order to formulate Proclus' argument accurately in the new constructive setting let us first formulate and discuss the following constructive version of \textbf{RU}:  

\begin{equation}
\begin{array}{l}
\text{Given objects $a, b$ and relation $R$ such that $aRb$ there exists \textbf{unique} object }\\ \text{(construction) $e(a,b,R)$ which makes it \emph{evident} that  $aRb$ holds, in symbols $e : aRb$}.
\label{eq4}
\end{array}
\end{equation}

Applying \eqref{eq4} to Proclus' argument we get the following reconstruction.  \eqref{eq4} tells us that there is at most one angle of form $\widehat{a, b}$ that can  witness the fact that $a\bowtie b$ (in symbols, $\widehat{a, b} : a\bowtie b$). But the cone construction demonstrates us that there exist two different angles of this form, $\widehat{a, b}_{P}$ and $\widehat{a, b}_{C}$, that can do this job. Hence an angle cannot be defined as (an instance of) relation $\bowtie$. 

Unlike the classical version of \textbf{RU}  \eqref{eq1} where the reference to relation $R$ turns out to be redundant, in the constructive version of the same principle \eqref{eq4} this reference is essential because evidence $e$ is proposition-specific and hence relation-specific\footnote{A point of intersection of two given straight lines can evidence that the two lines intersect but it cannot evidence that the lines are perpendicular.}. So \eqref{eq4}, unlike \eqref{eq1}, does not reduce to the simple claim \eqref{eq3} according to which two given objects always form a unique pair of objects. Thus \eqref{eq4} is certainly a more interesting principle than  \eqref{eq1} or  \eqref{eq3}. But is \eqref{eq4} really justified as a general principle? 

Apparently it is not. As any other proposition, a proposition of the form $aRb$ may have more than one proof. For example, two intersecting circles may have two different points of intersection, and each of these two points can serve one as an evidence that the given circles intersect. For a more up-to-date example consider the relation of isomorphism between algebraic groups. Two groups $G, G'$ are isomorphic, in symbols $G \sim G'$,  when there is an isomorphism between them, which is a structure-preserving map of form $i :G \xrightarrow{\sim}G'$. But such isomorphism $i$, if it exists, is not necessarily unique. Both these examples apparently refute \eqref{eq4}. It appears that \eqref{eq4} is plainly false, and hence Proclus' argument under the above constructive reconstruction is not conclusive.

Before I propose a tentative justification for \eqref{eq4} that will help me to save Proclus' argument in a constructive setting, let me make a terminological remark concerning the above example of isomorphic groups, which I'm going to reuse in what follows. In the modern colloquial mathematical language term ``isomorphism'' is used in these two closely related senses: it may stand (i) for a binary relation between mathematical structures and (ii) for a map between such structures, which can evidence isomorphism in sense (i). This ambiguity may make no harm in some other contexts but in the present context we need to distinguish between (i) and (ii) properly. For this end I shall call isomorphism in sense (i) \emph{iso-relation}, and an isomorphism in sense (ii) an \emph{iso-map}\footnote{Referring to iso-relations and iso-maps in what follows I mean by default group isomorphisms but my arguments involving this example equally apply to other types of mathematical structure.}.

In order to justify  \eqref{eq4} we need to read it ``internally'' rather than ``externally'' as we did this earlier. I shall demonstrate the internal reading of  \eqref{eq4} using, once again, the example of isomorphic groups. Judgement $i :G \sim G'$ allows one to recognise iso-map $i :G \xrightarrow{\sim}G'$ as a conclusive evidence for proposition $G \sim G'$ but it does not help one to distinguish between different evidences $i$ and $i'$ that may equally well qualify for that role. Earlier we simply took it for granted (i.e., assumed it ``externally'') that iso-maps $i$ and $i'$ were different \textemdash\ just as did Proclus when he assumed that angles $\widehat{a, b}_{P}$ and $\widehat{a, b}_{C}$ appearing in his cone construction\footnote{That is, angles $\angle AOB_{P}$ and $\angle AOB_{C}$ at Fig.\ref{fig:cone}} were different. But in the present constructive setting proposition $i\neq i'$ needs to be constructively justified (via a constructive refutation of $i = i'$) separately with a separate evidence. 

This can be done internally in an appropriate constructive framework (like MLTT) via a judgement of form $\alpha : i\neq i'$. But observe that $G \sim G'$  and $i\neq i'$ are two separate propositions, and $i :G \sim G'$ and $\alpha : i\neq i'$ are two separate judgements. Thus proofs $i,i'$ of the same proposition $G \sim G'$ cannot be distinguished without proving another proposition, namely  $i\neq i'$. Thus it is meaningful to postulate (as this is done in HoTT) that a single self-standing constructive proposition does not admit for different proofs but can only be either constructively true (in case it has a proof) or constructively false (in case that it provably has no proof). Formally, if $P$ is a proposition then $e : P$ and $e' : P$ imply $e =_{P} e'$. Thus the constructive propositions just described are bi-evaluated like classical propositions\footnote{The authors of the HoTT Book use for this constructive conception of proposition technical term ``mere proposition'' \cite[p.103, Def.3.3.1.]{UF:2013}. In \ref{HoTT} we shall see that the uniqueness of proof of a mere proposition (which serves as the defining property of this latter concept) also has in HoTT a geometrical (to read homotopical) motivation. \label{mp}}. The constructive version \eqref{eq4} of \textbf{RU} is a special case of the above general principle where propositions have special form $aRb$. 

How the above argument allows one to get around the obvious objection against \eqref{eq4} that a proposition of form $aRb$ (as any other mathematical proposition) may have more than one proof? Here is a way to do this. Recall that in the classical setting one wholly abstracts oneself from how a given proposition is known and whether it is known at all: one simply assumes that it has a truth-value. In the constructive setting we want a \emph{proof-relevant} concept of proposition: we count proposition $P$ as being true when it has some proof $e$ and count it as being false when it provably has no proof. Making judgement $e : P$ requires presenting proof $e$ in an explicit form. But in the constructive setting too we are making a (more limited) step towards abstracting propositions from their proofs. Take the example of parallels. It is common in the Euclid-style geometry to talk about a pair of parallels $a\parallel b$ having it in mind that such a thing is \emph{constructible} but without presenting any particular construction for it. The needed construction $e : a\parallel b$ can be always provided by a special request, and there is more than one such possible construction $e$ for a given pair of parallels $a\parallel b$. Thus we have here a room between the notion of being (potentially) constructible and the exhibiting of a particular construction here and now.

This is where the constructive notion of proposition explained above lives. Postulating (with \eqref{eq4}) counterintuitively that a given pair of parallels $a\parallel b$ admits for unique construction $e$ (that provides for judgement $e : a\parallel b$), we abstract away details of this construction that may allow us to distinguish $e$ from another construction $e' : a\parallel b$ that has exactly the same effect but can be otherwise very different. This abstraction allows us to think and reason of pair $a\parallel b$ as being \emph{constructible} without specifying any particular construction for it. This shows that the very idea of constructive mathematical reasoning admits for nuances and degrees: a stronger version of mathematical constructivism may rule out the abstract notion of potential constructibility (and the notion of ``construction by a special request'') on epistemological or other philosophical grounds, and require an explicit construction in all cases. Since the purpose of the present Chapter is to explain Proclus' argument in terms proposed by Vladimir in his message \ref{VR1}, I use here, somewhat dogmatically, the notion of being constructive as it appears in the standard version of HoTT \cite{UF:2013} without trying to explore other constructive approaches.

Let me now parallel Proclus' argument using the example of isomorphic groups. I shall argue that an iso-map can \underline{not} be properly defined as a pair of isomorphic  groups $G \sim G'$ provided with evidence $i$ that they are isomorphic (where $i$ is iso-map of the form $i :G \xrightarrow{\sim}G'$). Once again, one can 
say that I'm knocking here on an open door, and remark that the correct order of definitions is the opposite: one should first define the concept of iso-map (for groups) and only then introduce the concept of iso-relation by saying that two groups hold this relation when there exists an iso-map between them. This remark is fair enough but, once again, it misses the point of the present discussion. Let us assume for the sake of the argument that our constructive definition of the iso-relation for (a chosen class of) groups comprises an effective procedure (algorithm) $E$, which for any pair of groups $G, G'$ (from the chosen class) either constructs an iso-map $i_{E} :G \xrightarrow{\sim}G'$ between them or gives us an evidence that no such iso-map exists\footnote{
This strong assumption implies two things. First, it implies that the groups should, in their turn, be introduced constructively, say, with their finite presentations. Then an algorithm (program) that builds an iso-map between two given groups so presented will use syntactic details of their presentations. The ``choice'' of iso-map in this case will depend on syntactic choices made by the programmer. Second, recall that the group isomorphism problem in the general case of finitely presented groups is provably undecidable \cite{Stillwell:1982}. So in the general case no such universal procedure $E$ exists. But the problem is decidable in some special cases including the obvious case of finite groups presented with their full Cayley tables. For the sake of the example let us assume that we are talking here about the class of finite groups.}. Now the question is whether or not $E$ is sufficient for defining the concept of iso-map for groups (of the given class). One may plausibly argue in favour of the positive answer to this question by saying that since $E$ provides for the explicit construction of iso-map in all relevant cases (and moreover provides for an evidence of the non-existence of such maps in other relevant cases), it effectively introduces the concept of iso-map (for groups of the given class) along with that of iso-relation, so the introduction of \emph{term} ``iso-map'' becomes a matter of mere terminological convention. 

But this latter argument is erroneous. True, $E$ shows us how the generic iso-map looks like in the given case. This information is all what one needs to know about the defining shared properties of the iso-maps of the given class. What $E$ does not tell us, however, is how to identify equal iso-maps and distinguish between different iso-maps. Thus our constructive definition of iso-relation allows one to introduce the relevant general concept of iso-map in a sense but it does not support reasoning about iso-maps as identifiable mathematical objects\footnote{This peculiar situation makes echo of Immanuel Kant's systematic distinction between general concepts, on the one hand, and mathematically \emph{constructed} concepts, one the other hand. The ``construction of concepts'' (Konstruktion der Begriffe) grounded in the temporal and spatial mathematical intuitions is, in Kant's view, the main characteristic feature of mathematical reasoning which distinguishes it from a philosophical speculation, see \cite[A713/B741 ff.]{Kant:1999}. In Kantian terms our example can be analysed as follows. Procedure $E$ \emph{assumes} a construction of group concept (via its finite presentation) and effectively \emph{constructs} the concept of iso-relation for groups. This latter construction involves the \emph{concept} of iso-map, which is made explicit with $E$. At the same time, $E$ falls short of \emph{constructing} the concept of iso-map.} This is not acceptable for a constructive definition of a mathematical concept: a sound definition of this sort should also tell us how to identify and to distinguish between objects that fall under the defined concept\footnote{An obvious philosophical reference appropriate for supporting and analysing the claim is Willard Quine's popular dictum ``No entity without identity''  and his related considerations in \cite[Ch.1]{Quine:1969}. A combination of Kant's and Quine's insights brings about the following picture. Given a concept, one is always in a position to speculate about putative entities falling under this concept (unless the given concept is obviously self-contradictory). But in order to \emph{construct} the given concept, and use it in mathematics and in mathematically-laden science, one also needs to learn how to identify and distinguish such putative entities. While Kant believed that this can be always done with a spatio-temporal intuitive representation of a given concept (provided that the concept is so representable), Quine believed that the job can and should be done with a first-order theory of \emph{classes}, i.e., with some version of axiomatic Set theory \cite[p. 21]{Quine:1969}. HoTT provides a novel perspective on this traditional philosophical issue, which I leave here for a future research.} 
This is a reason why an iso-map can \underline{not}, after all, be properly defined as a pair of isomorphic objects $G, G'$ provided with procedure $E$ that constructs a witness for iso-relation $G \sim G'$ in form of iso-map $i_{E} :G \xrightarrow{\sim}G'$.

\emph{Mutatis mutandis} Proclus\footnote{modulo my proposed constructive reconstruction of his reasoning} applies the same argument to Euclid's Def.1.8., which introduces the concept of angle as (an instance of) binary relation $\bowtie$ of inclination between two given lines (i.e., the angle's sides). In the quoted fragment of his commentary on Euclid's Def.1.8. Proclus does not point to certain angles which do not instantiate this relation. Neither he points to certain instances of   $\bowtie$ which are not angles. His critical argument against Def.1.8. is compatible with the assumption according to which every instance of $\bowtie$ corresponds to an angle and every angle corresponds to an instance of $\bowtie$. The problem of Def.1.8 pointed to by Proclus is different: the aforementioned correspondence between instances of $\bowtie$ and angles is not ono-to-one but one-to-many, so a single instance of $\bowtie$ like $AO \bowtie BO$ at Fig.\ref{fig:cone} may correspond to two different angles $\angle AOB_{P}$ and $\angle AOB_{C}$. This feature of Def.1.8. disallows a direct identification of angles with instances of $\bowtie$. 

Proclus' argument is also compatible with the assumption according to which Def.1.8. correctly determines all essential features of the angle concept. At least, Proclus does not challenge this assumption in his quoted argument. What this definition fails to do, according to this argument, is to identify and to distinguish individual angles properly.

We can now see that the constructive interpretation of Relational Uniqueness given in the current Section provides for a subtler version of Proclus' argument than the classical interpretation presented earlier in \ref{RUclass}. According to the classical interpretation of this argument, an angle is not a pair of lines satisfying a certain condition simply because one and the same pair of lines satisfying this condition does not determine an angle uniquely. This classical interpretation leaves it unexplained where the requirement of uniqueness comes from. As it has been argued in the introductory part of \ref{RU}, a valid definition of a mathematical concept is not supposed to determine a unique object that falls under this concept. The constructive interpretation provides an answer to this question in the form of explicit general principle according to which a valid definition of a mathematical concept should allow one to correctly identify and to distinguish objects that fall under this concept. As shows Proclus, Euclid's Def.1.8. does not meet this requirement because the method of identification of angles provided with this definition is in odds with the common intuition that tells us that $\angle AOB_{P}$ and $\angle AOB_{C}$ at Fig.\ref{fig:cone} are two different angles \textemdash\ while Def.1.8. forces one to treat the two angles as the same.

 \subsection{Homotopy Type theory and Univalent Foundations of mathematics} \label{HoTT}  
The reader familiar with the basics of HoTT and the Univalent Foundations of mathematics (UF) \cite{UF:2013}, \cite{Grayson:2018} can already see that our constructive interpretation of Proclus' argument \ref{RUconstr} is inspired by HoTT and can be used for illustrating some semantic aspects of MLTT and HoTT. As shows are exchange with Vladimir \ref{AVC} I did not immediately understand what he meant by a distinction between propositions and structures in Proclus \ref{VR1}. But I got his point later after receiving his second message \label{VR2} where Vladimir replied to my direct request to explain what he meant by a mathematical structure \ref{RV1}. 


Let me now provide some hints for the reader not familiar with HoTT. The standard version of HoTT \cite{UF:2013} can be described as the constructive aka intuitionistic type theory due to Martin-L\"of (MLTT) \cite{Martin-Lof:1984} interpreted in terms of Homotopy theory , which is a part of Algebraic Topology that studies algebraic properties of continuous \emph{paths} in a topological space and their  \emph{homotopies}, i.e., ``paths between paths''. This ladder of paths, paths between paths, etc., can be continued upward indefinitely \cite{Gray:1975}. 

Unlike the standard Predicate Calculus (whether classical or intuitionistic) MLTT does not dispose of a notion of predicate but it comprises the notion of \emph{dependent type}, which plays a similar role. Given type $A$ one can have a family of types of form $B(x)$ where $x$ is a term of type $A$, in symbols

\begin{equation}
x : A \vdash B(x) : TYPE
\label{eq6}
\end{equation}
In particular, given type $A$ and two terms $t, t'$ of that type, one has the \emph{identity type} that depend on these two terms: If inhabited, the identity type shows that $t,t'$ are equal:

 \begin{equation}
t, t' : A \vdash t=_{A}t' : TYPE
\label{eq7}
\end{equation}      

The equality (aka identity) in \eqref{eq7} is called in MLTT/HoTT \emph{propositional} equality; at the first approximation type $t=_{A}t'$ can be thought of as a proposition and its terms can be thought of as proofs of this proposition. (We shall briefly see how HoTT extends this intended semantics of MLTT.) Proposition $t=_{A}t'$ along with its proof $p$ qualifies as a \emph{judgement} and has form $p : t=_{A}t'$\footnote{Notice that in MLTT and HoTT only terms of the same type can be compared, hence the presence of subscript $A$ in $t=_{A}t'$.}. The propositional equality should not be confused with a different kind of equality used in MLTT, which is called \emph{definitional} or \emph{judgemental}: as the name suggests judgemental equalities are judgements on their own and don't admit for or require a proof. At the first approximation the definitional equalities can be thought of as notational conventions. In what follows we discuss only propositional equalities. 

The construction of identity type in \eqref{eq7} can be iterated: given two terms $p,p'$ of type $ t=_{A}t' $, i.e., two proofs of (evidences for) proposition $ t=_{A}t'$, one can further compare the proofs:
 
 \begin{equation}
p,p' : t=_{A}t'  \vdash p=_{t=_{A}t'}p' : TYPE
\label{eq8}
\end{equation} 

The iteration of identity type can be continued indefinitely. Types like $p=_{t=_{A}t'}p'$ are referred to as \emph{higher} identity types. Homotopy theory helps to grasp the hierarchic syntactic structure of higher identity types in intuitive geometrical terms: think of type $A$ in \eqref{eq7}-\eqref{eq8} as a space, of terms $t,t'$ as points of this space, and of $p,p'$ in \ref{eq8} as paths between these points. The next iteration that compares these paths to each other involves ``paths between paths'', i.e., path homotopies. The above homotopical interpretation extends to all operations of MLLT; in particular the type dependence \eqref{eq6} is interpreted in HoTT in terms of \emph{fibration}, which is another basic concept of Homotopy theory.

It may happen that for given type/space $A$ the above iterative construction of (higher) identity types stops at certain step $n$ in the sense that at all further steps $k>n$ one always gets a unique proof that terms of $k$th identity type (obtained during the iterative process) are equal for all pairs of terms. 
In other words, it may happen that after $n$-th iteration the structure of (higher) identity types becomes trivial. This consideration provides us with a hierarchical classification of type/spaces by their \emph{homotopy levels}. If the iterative process stops at the very beginning, i.e., if $A$ is such that that for all $t,t'$, ${t=_{A}t'}$ is either empty, or has a single evidence, one qualifies ${t=_{A}t'}$ as a proposition (or more precisely a\emph{mere} proposition \footref{mp}, \cite[p.103, Def.3.3.1.]{UF:2013}, which is either true or false\footnote{Compare comments on \eqref{eq4} in \ref{RUconstr}.}. In this case the equality in \eqref{eq7} functions like the standard binary relation in the Predicate Calculus\footnote{In the sense that in both cases propositions are bi-evaluated, many other features still being different.}: given two terms $t,t'$ of some type $A$, relation $t=_{A}t'$ either holds or does not hold. Notice that the base type $A$ in this case can be described as a \emph{set}: terms of $A$ are either equal or not equal,  while all higher identity types beginning with $p=_{t=_{A}t'}p'$ (for all such $p,p'$) are all trivial (i.e., also propositional).     

If the iteration of identity types stops only at the next step then the situation changes dramatically. In that case type ${t=_{A}t'}$ may have two or more distinguishable proofs $p,p'$. Now terms $t,t'$ can be equal ``in different ways''  meaning that there can exist more than just one path that joins identical points $t,t'$.  Thus we get an example of what Vladimir in his message \ref{VR2} calls a \emph{structure}:  two elements $t,t'$ linked in two different ways $p,p'$ 
\footnote{In the Homotopy theory and in HoTT such a structure can be described as the \emph{fundamental groupoid} of the underlying space $A$. 
The assumption that $p=_{t=_{A}t'}p'$ and all higher identity types built on $A$ are trivial translates into the language of homotopies as follows. The identification of paths in $A$ (but not of points in $A$!) is a yes-no matter: paths count as equal iff they are homotopic, i.e., if there is a homotopy between them, or unequal iff there is no such homotopy. Here we do not further distinguish between different homotopies between the same pair of paths. At the next homotopy level one distinguishes between different path-homotopies but not between higher homotopies. \label{fn}}. 
As Vladimir specifies in the same message, relation is a special type of structure where the elements are linked in at most one way (i.e., either linked or not linked).

Notice that the homotopical hierarchy of types in HoTT concerns all types but not only identity types; the structure of identity types outlined above is a key for describing this general structure. In particular, it concerns those types that can be identified as ``binary relations''. Following Vladimir's suggestion, by relations in HoTT I understand  ``mere relations'', that is, families of dependent propositional types indexed by two variables (i.e., by terms of the base type) like in the case of propositional identity types of the form $t=_{A}t' $

The analogy with Proclus' argument presented above in \label{PE} now becomes obvious and straightforward. Proclus' cone construction shown at Fig.\ref{fig:cone} comprises two different angles $\angle AO, BO_{P}$ and $\angle AO, BO_{C}$, which serve as ``links'' between the two lines $AO, BO$ satisfying the conditions of  Euclid's Def. 1.8.. Each of these two angles witnesses the fact that lines $AO, BO$ hold the relation of \emph{inclination} $\bowtie$ referred to in this definition. Thus we have here an instance of \emph{structure} in Vladimir's sense: two objects $AO$ and $BO$ are linked in two different ways. This structure is not a mere relation because objects $AO$ and $BO$ are linked in more than one way. This HoTT-inspired interpretation of Proclus' argument extends to its conclusion according to which angle is not a (mere) relation.

 


 \subsection{Philosophical and Historical Remarks, and Some Pointers to a Future Work} \label{FW}
 \subsubsection{Mathematical Structuralism and the Set-Level Mathematics} \label{SLM}
How the conception of structure presented in Vladimir's message \ref{VR1} compares to the standard (modulo its multiple variations) notion of mathematical structure as a family of sets with some relations \ref{RSV}? Here is a technical part of the answer. Recall that in the case when type $p=_{t=_{A}t'}p'$ is a mere proposition (i.e., has at most one term) the base type $A$ is a \emph{set} (in the sense of HoTT)\footref{fn}. Now if the family of types $B(x)$ dependent on type $A$ (as in \eqref{eq6}) is such that for all $x$, $B(x)$ is also a proposition then $B(x)$ can be thought of as a monadic predicate aka a \emph{property} that defines a subset $B$ of set $A$ \cite[p. 106-107]{UF:2013}. In a similar way HoTT allows for mimicking binary predicates and predicates with larger arities, i.e., relations, which are defined on families of sets and thus form analogues of the usual set-based mathematical structures. Such analogues of the set-theoretic structures feature in HoTT is a very special case of general HoTT-structures referred to by Vladimir in our exchange: generally, the base type $A$ is not necessarily a set (but can be a type of a higher homotopy level) and dependent types like  $B(x)$ (where $x : A$) are not necessarily (mere) propositions. 
 
Thus the notion of structure in HoTT is by far more general than the standard notion. Vladimir's personal motivation in favour of this generalisation, as I understand it, was to bypass the limits of the standard ``set-level'' mathematics and develop a new kind of mathematics beyond this limit. How this new mathematics could look like Vladimir knew from his work in the Algebraic Geometry that won him the Fields Medal in 2002: his whole project of designing new ``univalent'' foundations of mathematics on the basis of HoTT aimed at a novel formal framework for mathematics that could make proofs in Algebraic Geometry and other highly technical fields of today's abstract mathematics more rigorous, more transparent, more accessible for colleagues working in other mathematical disciplines, and verifiable with a computer \cite{Voevodsky:2014}. 
 
Vladimir was not alone who reflected on the notion of mathematical structure in HoTT back in 2016 and during the later years, see \cite{Awodey:2014},\cite{Chen:2024},\cite{Corfield:2017},\cite{Shulman:2017},\cite{Teruji:2014}, \cite{Tsementzis:2017}. A systematic review of this literature is out of place here but it is appropriate to say few words on how Vladimir's conception of mathematical structure presented in this Chapter compares to what other people working in HoTT say about the same subject. First of all, a reservation is here in order. Unlike some other working mathematicians having a serious interest to the history and philosophy of mathematics Vladimir was never engaged in academic discussions and publications in these fields. I doubt that Vladimir ever followed discussions on the Mathematical Structuralism published in the philosophical academic journals. My exchange with Vladimir on philosophical and historical matters is not an exception: it should not be read as Vladimir's systematic account of the concept of mathematical structure. If Vladimir were asked to write an encyclopaedia article on this subject he would certainly say more than he did in our brief exchange.

Having this reservation in my mind I nevertheless would like to share one observation. The continuing philosophical discussion on the concept of mathematical structure and Mathematical Structuralism in HoTT started with Steve Awodey's paper \cite{Awodey:2014} where the author shows how HoTT  makes true the old structuralist dream by providing a formal framework where isomorphic set-based mathematical structures can be rigorously (rather than only informally as it has been done earlier \cite{Bourbaki:1950},\cite{MacLane:1996}) called and treated as equal.  Awodey calls the statement ``Isomorphic objects are identical'' the ``Principle of Structuralism (PS)'' \cite[p.2]{Awodey:2014} and shows that it is implied by the Axiom of Univalence (AU)\footnote{PS is an informal linguistic description of the propositional version of AU, i.e., of AU applied to propositional types. A more detailed analysis of this feature of HoTT is given in \cite{Ahrens&North:2019}.} This allows him to qualify \emph{all} mathematical objects in UF (and in the UF-based mathematics) as structures and claim that this is the strongest possible formulation of Mathematical Structuralism. The same feature of UF serves as a starting point for discussing Mathematical Structuralism for David Corfield \cite{Corfield:2017},\cite{Corfield:2020}, Dimitris Tsementzis \cite{Tsementzis:2017} and other researchers.

Vladimir, of course, appreciated the fact that UF helped to make the conventional informal mathematical talk about the ``equality up to isomorphism'' formally justified and logically rigorous. Moreover, Vladimir's talk on the Foundations of Mathematics and Homotopy Theory given in the Institute of Advanced Studies (IAS) in Princeton on March 22, 2006 several years before he formulated the Axiom of Univalence and coined the title of Univalent Foundations in 2010 \footnote{Publicly this name first appeared in his talks in Bonn and in the IAS in 2010. Materials of Vladimir's talks mentioned in this Chapter are publicly available via his memorial page on the website of IAS at \url{https://www.math.ias.edu/vladimir/Lectures} .}
makes it clear that some form of PS \textemdash\ or, more precisely, the general problem of identity aka equality in mathematics where the PS belongs  \textemdash\  was one of Vladimir's motivations behind this project. But in his later public presentations of UF Vladimir rarely touched upon this issue focusing on different aspects of UF. Perhaps Vladimir believed that once the issue of identification of isomorphic structures was settled in HoTT with the AU, it did not deserve much further discussion.  In 2016 Vladimir gave a very interesting talk on Multiple Concepts of Equality in the New Foundations of Mathematics (Bielefeld, July 18, 2016) where he discussed the judgemental and propositional equalities in MLTT and in HoTT but, once again, he did not mention in this talk anything like PS or anything else that resembled the usual Structuralist agenda. Thus the available evidence suggests that this agenda played a very limited role in Vladimir's thinking about UF. It also explains why Vladimir does not refer to anything like PS talking about his general conception of mathematical structure in \ref{VR2}.

The technical feature of HoTT/UF  which Vladimir used in his definition of mathematical structure, belongs to the foundation of this theory \ref{HoTT}; every type in UF is a structure in Vladimir's sense of the term just like it is a structure in Awodey's sense of the term. So in the UF the two conceptions coincide in their extensions 
\footnote{Notice, however, that the conception of structure highlighted by Vladimir, unlike the standard conception stressed by Awodey and others, does not require the presence of AU.} 
But the two conceptions of structure differ nevertheless in their conceptual contents. While the conception of mathematical structure developed by Awodey and others in the novel framework of UF is obviously rooted in the Bourbaki-style Structuralism of the 20th century \cite{Bourbaki:1950}, the conception of mathematical structure put forward by Vladimir in the same context highlights a feature of HoTT/UF that hardly has any obvious counterpart in the Bourbaki-style mathematics.  

Concluding his discussion on the concept of mathematical structure in HoTT, David Corfield formulates the following methodological principle: 

\begin{quote}        
 Any time we have a construction which traditionally has been taken to apply only to sets or only to propositions, then since in HoTT these form just a certain kind of type, we should look to see whether the construction makes sense for all types. \cite[p.699]{Corfield:2017}        
\end{quote}

Corfield talks here about possible extensions of the familiar set-level structures over types of higher homotopy levels. This strategy comes close to Vladimir's motivation behind his notion of structure, as I understand it. But there is also an important difference. I believe that Vladimir's strategy of exploring the world of higher homotopy levels, which was based on his experience in the Algebraic Geometry, was rather this: to develop higher-level structures to begin with and then look back at the set-level mathematics as a source of elementary examples and toy models. Compare what Vladimir says about a ``laughably simplified version'' of a mathematical problem in his first message \ref{VR1}. Thus Vladimir could agree with Corfiled at this point but at the same time qualify his proposed research strategy as too careful and too conservative.

Let me now provide some critical comments on the idea of ``going beyond the set-level mathematics'' as it appears in the context of the UF. It should not be confused with the well-known metamathematical issues concerning the expressive and the proof-theoretic power of ZFC or another similar theory of sets. First, the concept of set as it features in HoTT and in the ZFC are not exactly the same concepts. Second, the limits of the set-level mathematics in the UF should be understood in terms of actual mathematical practice (both in the mathematical research and the mathematical education) rather than in abstract meta-theoretical terms as this is so often done in philosophical analyses of Set theory as a foundation. 

In term of UF, a non-trivial groupoid  is a structure of a higher homotopic level than the level of set-based structures. But such a groupoid can be also construed by the standard set-theoretic means as a (non-directed) multigraph with a partial operation of composition defined on its arrows. The fact that a groupoid can be so construed does not mean, however, that it should be. Instead of using the language of sets throughout the mathematics one can use the language of groupoids (having it in mind that sets are trivial groupoids), which brings about a different kind of mathematics. Such mathematics is different at least in the sense that its basic concepts are different and its theories are built in a different way.  Whether the obtained theory is translatable into the language of sets \emph{in principle} or not is quite a different issue. In many cases of interest such a translation is possible but this meta-theoretical fact does not make UF less important or less interesting for a philosophical analysis. Similar remarks can be made about all higher groupoids up to the  $\infty$-groupoids.  

Seen from the new perspective opened by the UF, the choice of the concept of set as \emph{the} elementary (primitive) concept for all mathematics appears somewhat arbitrary; it appears that this choice can be explained by metaphysical and aesthetic arguments rather than by scientific and properly mathematical ones. It goes without saying that Set theory played an important role in the 20th century mathematics as this is evidenced, in particular, by David Hilbert's famous ``Cantor's Paradise'' metaphor \cite[p.170]{Hilbert:1926}. But the UF provides us today with a new foundational perspective, which gives room for a more critical view on the role of Set theory in the history of the 20th century mathematics. 

The idea of developing new mathematics beyond the limits of the set-based mathematics is older than UF by at least sixty years. As it has been already mentioned in \ref{RSV}, Bourbaki, who were probably the most important promotors of the set-based mathematics in the 20th century, at the same time urged in their 1950 manifesto \cite{Bourbaki:1950} for the ``disappearance of sets'' in favour of abstract mathematical structures. During the second half of the 20th century this idea was associated with the practice of using the ``language of categories'' (in the sense of the mathematical Category theory) instead of (or along with) the more familiar ``language of sets'' in many advanced mathematical disciplines including Homological Algebra, Algebraic Topology and Functional Analysis \cite{Manin:2002}. In the mid-1960s William Lawvere started a project of developing new foundations of mathematics on the basis of general Category theory \cite{Lawvere:1964},\cite{Lawvere:1966} but since its critics argued that Category theory could not be developed without using some concept of set or collection at the foundational level \cite{Feferman:1977}, the issue always remained controversial. 

Vladimir's first attempts to develop new foundations of mathematics allowing to reach ``beyond'' the set-based mathematics were also motivated by the general Category theory. But in his 2014 lecture where Vladimir tells his personal story of inventing the Univalent Foundations, he, surprisingly, describes this popular idea as a ``roadblock'':
  
\begin{quote}   
The greatest roadblock for me was the idea that categories are ``sets in the next dimension''. I clearly recall the feeling of a breakthrough that I experienced when I understood that this idea is wrong. Categories are not ``sets in the next dimension''. They are ``partially ordered sets in the next dimension'' and ``sets in the next dimension'' are groupoids.

This new perspective on ``groupoids'' and ``categories'' took some adjustment for me because I remember it being emphasised by people I learned mathematics from that one of the things that made [Alexandre] Grothendieck's approach to algebraic geometry so successful was that he broke with the old-schoolers and insisted on the importance of considering all morphisms and not only isomorphisms. (Groupoids are often made of set-level objects and their isomorphisms, while categories are often made of set-level objects and all morphisms.) \cite{Voevodsky:2014}
\end{quote}

Indeed, the idea of extending the familiar set-level mathematics to the groupoid-level mathematics to the mathematics of higher-groupoid levels is proper to HoTT/UF. General categories are construed in the UF-based mathematics on the top of groupoids along with other types of mathematical structure. I must confess, however, that I'm not convinced along with Vladimir (back in 2014) that the choice between general categories and groupoids as elementary mathematical structures has been finally settled with the development of UF in its standard form. There are strong conceptual reasons to believe that  Alexandre Grothendieck might be right, after all, that general morphisms in a category should be treated on equal footing with isomorphisms \cite{Rodin:2011} (so the ``old-schoolers'' referred to by Vladimir in the above quote might still be wrong). The recent work on the \emph{directed} version of HoTT (Directed Type theory or DTT for short) gives a reasonable hope to build a new type-theoretic formal framework with semantics in general higher categories rather than specifically in higher groupoids like in the standard HoTT. By the date, this project has been already realised for some types of higher categories\footnote{For a general informal introduction to DTT, its philosophical significance, and further references to mathematical sources see \cite{Rodin:2024}.}.

  \subsubsection{Frege and Lotze} \label{FL}
Back in 2016 I did not take very seriously Vladimir's historical question about the date and the place when mathematicians and philosophers, as Vladimir suggested, privileged relations over structures in their basic world-picture \ref{VR1}. Now I can see that I was wrong. Let me explain this. 

Since Vladimir's reading of Proclus was overtly anachronistic his historical question did not appear to me well-posed. I tried to answer this question by pointing to the fact that thinking of relations as non-monadic predicates in the sense of today's First-Order logic is, historically speaking, a relatively recent idea that goes back to the pioneering works by Gottlob Frege in the second half of the 19th century but hardly further back into the history. Since Vladimir's question also concerned ontology and metaphysics \textemdash\ or so I understood his words about the ``foundation of the world'' in  \ref{VR1} \textemdash\ I recalled of Betrand Russell who, according to his own testimony, designed a metaphysics on the basis of the Frege-style logic rather than the other way round
\footnote{Russell:
\begin{quote}
As I have attempted to prove in \emph{The Principles of Mathematics}, when we analyse mathematics we bring it all back to logic. $\dots$ In the present lectures, I shall try to set forth in a sort of outline, rather briefly and rather unsatisfactorily, a kind of logical doctrine which seems to me to result from the philosophy of mathematics \textemdash\ not exactly logically, but as what emerges as one reflects: a certain kind of logical doctrine, and on the basis of this a certain kind of metaphysics \cite[p.495-496]{Russell:1918}  
\end{quote}
}. 

In his reply \ref{VR2} Vladimir disagreed with me that the alleged oblivion of structures and the preference of relations in the common contemporary world picture could come from logic, and he pointed to \emph{Metaphysics} by Hermann Lotze \cite{Lotze:1841}\footnote{I didn't check this with Vladimir but I believe that he could read Lotze's Metaphysics in English translation \cite{Lotze:1884}} where Vladimir found the idea of ontology based on objects and relations described informally in traditional philosophical terms. In the same message Vladimir also remarked  that he found similar ideas in George Boole. I was surprised to see Lotze on Vladimir's reading list and looked only very superficially into the recommended texts, so this historical discussion did not then bring us anywhere because of my laziness.  

Today, after a reflection and some more reading, I can see that Vladimir was perfectly right. Back in 2016 I rather uncritically shared the popular interpretation of Frege's work in the line of Michael Dummett \cite{Dummett:1973} who believed that tracing the historical genesis of Frege's thought was not an appropriate way of its understanding. I was not  aware then about the controversy between Dummett and Hans Sluga concerning the relevance of historical questions in reading Frege. Today I find myself on Sluga's side of this debate. A more recent scholarship on Frege convinced me that Lotze's philosophy is crucial for understanding where Frege's logical ideas come from \cite{Gabriel:2002},\cite{Heis:2013}. 

I still have no idea how Vladimir picked up on Lotze. Some of Vladimir's remarks that I have not included in the present Chapter suggest that he was not sure about dates and didn't know that Frege was Lotze's student. In his message \ref{VR2} Vladimir wrongly suggested that Lotze did not read Boole. Anyway Vladimir somehow identified historical sources that were most relevant to his historical query. This convinced me later that Vladimir's historical intuition was in fact as strong and as sharp as his mathematical intuition. 

Thus my replies \ref{RV1} and \ref{RV2} to Vladimir's historical questions I judge today to be rather shallow; at best they can serve as a starting point of a more serious research on the historical genesis of the modern concept of relation in the new theoretical perspective provided by the Homotopy Type theory. The very fact that Vladimir's anachronistic reading of Proclus is coherent is quite remarkable: it shows that Proclus' Aristotelian conception of relation and the modern Fregean conception of relation as a predicate share an important common feature, which would be difficult or probably even impossible to identify without contrasting relations to structures as they appear in HoTT. As long as the concept of mathematical structure is understood in the non-standard way suggested by Vladimir, the question about its long-term history makes perfect sense. The existing historical accounts of the rise of Mathematical Structuralism \cite{Reck&Schiemer:2020} provide only a partial answer to this question. A lot more remains to be done to give to Vladimir's historical question a better answer. I leave this issue for a future research.

\section{Fragments of Correspondence between Andrei Rodin and Vladimir Voevodsky (January-February 2016)} \label{AVC}

The following are fragments of four electronic messages that Vladimir and I exchanged in the beginning of the year 2016. I do not quote the remaining parts of these messages because they include some unrelated discussions, and I opt for focusing this Chapter on a single topic. In order to facilitate the reading of these messages I include some references and some comments as footnotes. The original language of our correspondence was Russian; the English translation is mine. The words in square brackets [] are added to Vladimir's and my own words in order to clarify their meaning.

 \subsection{Voevodsky to Rodin, January 27, 2016}\label{VR1}

[$\dots$] In Proclus we find a clear distinction between a property and a structure (see his reflection on how to define an angle in his Commentary on Euclid)\footnote{See \emph{Commentary to the First Book of Euclid's Elements}  \cite{Proclus:1873}, English translation \cite{Proclus:1970}. Vladimir refers here to Proclus' commentary on Euclid's Definition 1.8. of plane angle, see  \cite[p.153]{Euclid:1908}, found in  \cite[p.98-104]{Proclus:1970}. More specifically, Vladimir refers to the fragment of this commentary \cite[p.99]{Proclus:1970} quoted above in \ref{PE}.}. Relation is a joint property of two or more objects. In philosophy and in mathematical logic the model of  the ``world'' is based on a collection of objects and a collection of relations between these objects, i.e., a collection of properties of assemblies of objects. I wonder when and how there emerged the strange idea that the foundation of the word can be described with a collection of relations between objects rather than with a collection of joint structures on assemblies of objects. There should be a [historical] moment when this idea was first presented as an auxiliary simplification : [as if someone says] ``let's assume this laughably simplified version''  just like in mathematics people often consider a simplified version of a problem in order to test with it this or that general idea. 

Volodya

P.S.  A relation either holds or doesn't hold (straight lines are either parallel or not) but a structure can have more than one representative. 

 \subsection{Rodin to Voevodsky Janaury 30, 2016}\label{RV1}

[$\dots$] [Replying to Vladimir's sentence ``There should be a [historical] moment when this idea was first presented as an auxiliary simplification.'' in \ref{VR1}].

In my understanding, this idea first emerged in logic (in Frege) and later Russell designed for this logic an appropriate ontology. At least this is how Russell himself describes this story in the Introduction to his [Philosophy of] Logical Atomism of 1918\footnote{See  \cite[p.495-496]{Russell:1918} and the discussion in \ref{FL}.} 
But in the context of the contemporary formal and mathematical logic this [move] did not look like a simplification because there was no any other formal theory [i.e., no alternative logical calculus] for working with relations. There were, however, interesting informal discussions, including the argument between Russell and F.H. Bradley on the ``internal and external'' relations. In his \emph{Logic} \cite{Bradley:1922} Bradley, using some Hegel's ideas, develops a very different theory of relations claiming that relata of a given relation in some (rather obscure) sense make with the relation one whole. Russell believed that this was an unnecessary complication. His main argument, as I understand it, was that his theory of relations was formalised while Bradley’s theory was not
\footnote{Of course, this is a polemical remark in passim but not an attempt on my part to give an account of or summarise the Russell-Bradley debate. A historical truth behind my 2015 remark is that Russell promoted using symbolic mathematical methods in the philosophical logic, while Bradley opposed it. I believe that the role of mathematical methods in the philosophical logic \textemdash\ and, conversely,  the role of philosophical reflection in the mathematical logic \textemdash\ need more methodological reflections than it is presently given. My current position is that mathematical methods are indispensable for fixing logical principles and ideas but a traditional form of philosophical reflexion on logical matters using the natural language is equally indispensable for building and shaping formal logical calculi and giving them a meaning. It goes without saying that we find in Russell a lot of informal philosophical discussions concerning logical formalisms. But my concern is that such discussions often miss a critical aspect: the author aims at providing a new logical calculus with a firm philosophical grounding rather than to discovering its epistemological limits and conceiving of possible alternatives. For a detailed analysis of the Russell-Bradley debate and of different traditions of its interpretation see \cite{Candlish:2007}. See \ref{FL} for discussion.}. This is hardly an essential objection. I doubt that Bradley’s theory [of relations] has anything to do with what you call a structure but this would be interesting to check. 

Could you clarify what you call a structure? If it is not a set with a collection of relations then what it is? [$\dots$]    A structure can be also understood as a set with relations [identified] up to isomorphism  \textemdash\ in this case a structure can have more than one representative [that is, more than one instantiation]. Does this concept of structure fit to what you’re talking about\footnote{Thus I first wholly misunderstood Vladimir's reference to Proclus in \ref{VR1} and thought of ``structures'' referred to by Vladimir in the usual way as ``mathematical objects identified up to isomorphism'' like ``the'' infinite cyclic group $\mathbf{Z}$. This is why I suggested that by different ``representatives'' of a given structure Vladimir might mean different isomorphic copies of the same mathematical structure .}?

Andrei 

 \subsection{Voevodsky to Rodin, January 31, 2016}\label{VR2}

In the context of my message [\ref{RV1}] a structure on several objects is an entity that can link these objects in a number of different ways.  [Given a pair of objects for which a given relation is well-defined] the relation [between them] either holds or does not hold, i.e., either there is a link [between the two objects] or not, so the set of possible versions of this link is either empty or has one element. A structure is a link such that the set of its possible versions can have more than one element. Notice that given this definition [of structure], a relation between $A$ and $B$ is a special case of structure on collection \{A,B\}. 

I doubt that the idea to use relations [instead of general structures] first appeared in logic. Look at [Hermann] Lotze's Metaphysics \cite{Lotze:1841} [English translation \cite{Lotze:1884}]: in the first 100 pages of this book the author discusses a world picture based on a set of objects and their relations. [$\dots$] . By the way, in the first [George] Boole's work [we also find] objects and relations. Perhaps Boole provides a reference pointing to his source [of the idea of reasoning in terms of objects and relations]
\footnote{Unfortunately I didn't check it with Vladimir which Boole's work he referred to. I believe that Vladimir refers here to Boole's first published logical essay of 1847 \cite{Boole:1847} (see also the modern commented edition  \cite{Boole:2009}) rather than to one of Boole's mathematical papers, which have been published earlier. In this work  Boole uses term ``relation'' referring to \emph{equations} of his logical calculus. Since an equation is a proposition, this meaning of the term falls under Vladimir's conception of relation explained above in \ref{RSV}. It is worth noticing, however, that Boole's equational propositions that he calls relations belong to his metalanguage rather than to the object-language formalised with his logical calculus. For an account of historical genesis of Boole's logical works see \cite{Laita:1979},\cite{Laita:2005}.}. 
But I doubt that Lotze read Frege or Boole\footnote{Gotlob Frege was Hermann Lotze's student in Jena; there is firm historical evidence that Frege carefully read Lotze's \emph{Logic} \cite{Lotze:1843}, \cite[p.122]{Heis:2013}. I could not find an evidence in the current historical literature that Lotze, reciprocally, read Frege's works \textemdash\  or at lest read the review of Frege's \emph{Conceptual Notation} of 1897 published by Ernst Schr\"oder in 1881 \cite{Schroeder:1881}. In any event, Lotze's \emph{Metaphysics} published in 1841 \cite{Lotze:1841} and his \emph{Logic}  \cite{Lotze:1843} published in 1843 could not be possibly influenced by ideas and writings of his student Frege who obtained his first doctoral degree in Jena only in 1873. Both Frege and Lotze read and criticised Boole's logical works \cite[p.126]{Heis:2013} (so Vladimir wrongly suggests in this message that Lotze was unaware about Boole's work). Vladimir was apparently not aware about these historical details when he referred to Lotze and to Boole in this message. See \ref{FL} for discussion concerning the significance of Vladimir's historical remarks.}.

 Volodya

 \subsection{Rodin to Voevodsky, February 1, 2017}\label{RV2}

Now I understand what you are talking about, and see a relevance of Univalent Foundations, of course.[$\dots$] You are right that the idea of world as a set with relations is older [than the Predicate Logic]. [$\dots$]  But in my view, for answering the question ``Why the simplified world picture [as a set of objects with relations] is still popular today in the mathematical logic and the analytic metaphysics?'' Russell is particularly important.  [$\dots$] Since doing metaphysics without using any formal machinery (as Lotze and Bradley did this in the 19th century) in the logical community of the 20th century became unfashionable, this [Russell's] simplified metaphysics [that construes the universe in terms of primitive objects and their relations] became a dogma. In Russell \cite{Russell:1918} we can see the exact moment when this happened
\footnote{See the quote from \cite{Russell:1918} in \ref{FL}}. 

I am interested in the Univalent Foundations also because they allow us to revise this metaphysics without giving up the idea of working formally. Such a revision is necessary, among other things, because the construction of the world as a set with relations is [apparently] not suitable for representing what physics tells us about the world. Attempts to apply the standard logical methods in physics, which have been made in the 20th century, did not bring significant results. [But there is a hope that novel logical methods related to the UF may perform in physics and other sciences better.\footnote{The argument is developed in my \cite[Chs. 2.3, 4.2, 4.3]{Rodin:2020}.}]

Andrei

\bibliographystyle{plain} 
\bibliography{voepro} 

\end{document}